# A VICINAL SURFACE MODEL FOR EPITAXIAL GROWTH WITH LOGARITHMIC FREE ENERGY

Yuan Gao[*]

Department of Mathematics
Hong Kong University of Science and Technology
Clear Water Bay, Hong Kong

Hangjie Ji[‡], Jian-Guo Liu[#,†] and Thomas P. Witelski[#]

Department of Mathematics[‡]
University of California, Los Angeles
Los Angeles, CA 90095, USA

Department of Mathematics[#]
Department of Physics[†]
Duke University
Durham, NC 27708, USA

(Communicated by Shouhong Wang)

Abstract. We study a continuum model for solid films that arises from the modeling of one-dimensional step flows on a vicinal surface in the attachment-detachment-limited regime. The resulting nonlinear partial differential equation, $u_t = -u^2(u^3 + \alpha u)_{hhhh}$, gives the evolution for the surface slope $u$ as a function of the local height $h$ in a monotone step train. Subject to periodic boundary conditions and positive initial conditions, we prove the existence, uniqueness and positivity of global strong solutions to this PDE using two Lyapunov energy functions. The long time behavior of $u$ converging to a constant that only depends on the initial data is also investigated both analytically and numerically.

1. **Introduction.** Epitaxial growth of crystal surfaces below the roughing temperature has attracted extensive interest. Unlike traditional modeling for fluid and solid mechanics where one starts with continuum theories for macroscopic variables, the modeling of crystal films was first established from the atomic perspective. At the nanoscale, crystal surfaces consist of basic structures such as interacting line defects (steps) and flat surface regions (facets). With adatoms detaching from one step and reattaching to another step after traveling along the facets, a step flow describes the mass transport along the crystal surface. For broader physical surveys of crystal growth we refer readers to [3, 10, 21].

From both discrete and continuum viewpoints, different models have been constructed to characterize step flows. From the discrete perspective, the dynamics of the steps are described by the step velocities using a Burton-Cabrera-Frank (BCF) type framework [3], which was investigated by Ozdemir and Zangwill [20]. The motion of these individual steps are usually modeled by systems of differential equations for step locations. At the macroscopic scale the continuum limit of step evolution is generally reformulated as nonlinear PDEs [5, 9, 12, 17–20, 22, 26]. Based on the conservation of mass, the dynamic equation for the surface height of a solid film, $h(t, x)$, can be written as

$$\frac{\partial h}{\partial t} = \nabla \cdot \left( M(\nabla h) \nabla \left( \frac{\delta G}{\delta h} \right) \right), \tag{1.1}$$

where the mobility function $M(\nabla h)$ is a functional of the gradient of $h$ and $G(h)$ represents a surface energy.

The mobility function takes distinctive forms in different limiting regimes. In the diffusion-limited (DL) regime, where the dynamics is dominated by the diffusion across the terraces, $M$ is a constant, $M \equiv 1$. While in the attachment-detachment-limited (ADL) case, the dominant processes are the attachment and detachment of atoms at step edges, and the mobility function takes the form [1]

$$M(\nabla h) = |\nabla h|^{-1}. \tag{1.2}$$

One form of the surface energy associated with the dynamical model (1.1) was identified [4] as

$$G(h) = \int \alpha |\nabla h| + \tfrac{1}{2}|\nabla h|^3 \, dx, \tag{1.3}$$

where the first term represents the step-formation energy and leads to a conical singularity at $\nabla h = 0$, and the second term is the step-interaction energy and is consistent with the discrete model for the crystal surface which will be discussed later. In the DL regime, Giga and Kohn [14] rigorously showed that with periodic boundary conditions on $h$, finite-time flattening to a spatially-uniform solution, $h \equiv C$, occurs for $\alpha \neq 0$. A heuristic argument provided by Kohn [11] indicates that the flattening dynamics is linear in time. However, in the ADL regime with the nonlinear mobility given by (1.2), the dynamics of the surface height equation (1.1) is still an open question [11].

In the ADL case, using (1.2) and (1.3), the evolution equation for the surface dynamics becomes

$$h_t = -\nabla \cdot \left\{ \frac{1}{|\nabla h|} \nabla \left[ \nabla \cdot \left( \alpha \frac{\nabla h}{|\nabla h|} + \tfrac{3}{2}|\nabla h|\nabla h \right) \right] \right\}. \tag{1.4}$$

To simplify the problem, consider a one-dimensional monotone vicinal surface corresponding to a monotone step train in the discrete model. Without loss of generality, assume that $h_x > 0$ for the whole domain. Under this assumption the $\alpha$ term drops out and (1.4) can be rewritten as

$$h_t = -\left[ \frac{1}{h_x} \left( \tfrac{3}{2} h_x^2 \right)_{xx} \right]_x. \tag{1.5}$$

This equation can also be derived from the dynamic equation (1.1) using the surface energy (1.3) with $\alpha = 0$. Following the method introduced by Ozdemir and Zangwill [20] and Shehadeh et al. [1], the new variable $u(t, h) = h_x(t, x) > 0$ can be adopted to rewrite the PDE as an evolution equation for the surface slope,

$$u_t = -u^2 (u^3)_{hhhh}. \tag{1.6}$$

This is a fourth-order degenerate parabolic PDE and the existence of global weak solutions for this problem has been shown previously in [8]. Specifically, it was proved that the solution $u$ is positive almost everywhere and that $u$ converges to a constant solution in the limit of long times. The central difficulty is how to deal with singularities if $u$ touches down to zero.[1] A method for degenerate parabolic equations used in [2] was adapted by neglecting a zero measure set $P_T = \{h \in [0,1]; u(h) = 0\}$, which is a closed set so that we can define a distribution on $(0,T) \times (0,1) \backslash P_T$ and formulate the definition of weak solutions to (1.6). Using this definition, the existence and almost-everywhere positivity of global weak solutions of (1.6) was proved.

Note that the slope equation (1.6) can also be derived from discrete models in the BCF framework [3] by showing the convergence of the discrete model to its continuum limit [7]. Let $\{x_i(t), i \in \mathbb{Z}\}$ be the step locations at time $t$ and let the height of each step be a constant $a = H/N$, where $x_{i+N}(t) - x_i(t) = L$ is specified for a train of steps with a length scale $L$, and $H$ represents the total height of $N$ steps. In the ADL regime, the BCF model without deposition flux is expressed by the step-flow differential equations,

$$\frac{dx_i}{dt} = \frac{1}{a^2}(\mu_{i+1} - 2\mu_i + \mu_{i-1}) \qquad \text{for } i = 1, \cdots, N,$$

where $\mu_i$ is a chemical potential giving the gradient of the free energy with respect to changes in the steps, $\mu_i = \partial \mathcal{E}/\partial x_i$. When the free energy involves only local contributions due to the interaction among steps, $f(r)$, it can be written as

$$\mathcal{E} = a \sum_{i=0}^{N-1} f\left(\frac{x_{i+1} - x_i}{a}\right).$$

Then, defining the step slopes $u_i(t) = a/(x_{i+1}(t) - x_i(t))$, we obtain the differential equations for slopes $u_i$ for $i = 1, 2, \cdots, N-1$,

$$\frac{du_i}{dt} = \frac{1}{a^4} u_i^2 \left[ f'\left(\frac{1}{u_{i-2}}\right) - 4f'\left(\frac{1}{u_{i-1}}\right) + 6f'\left(\frac{1}{u_i}\right) - 4f'\left(\frac{1}{u_{i+1}}\right) + f'\left(\frac{1}{u_{i+2}}\right) \right]. \tag{1.7}$$

Note that the right hand side of (1.7) is equivalent to a centered finite difference discretization in the limit of the step height $a \to 0$, or equivalently, as the number of steps $N \to \infty$. Therefore the solution of the slope ODE (1.7) should converge to the solution $u(t, h)$ of the continuum model

$$u_t = u^2 \left[ f'\left(\frac{1}{u}\right) \right]_{hhhh}, \tag{1.8}$$

where the step slope $u$ is considered as a function of $h$. For vicinal surface models with entropic and elastic-dipole interactions, the continuum model (1.8) is equivalent to the PDE (1.6) in [8] with the local contribution $f$ being $f_1(r) = \frac{1}{2}r^{-2}$. In 2002, Xiang [26] investigated a continuum model in the DL regime with nonlocal contributions with $f_2(r) = -\ln|r|$.

---

[1] In [8] the solution $u(t,h)$ of (1.6) was regularized by a small $O(\epsilon)$ perturbation, and it has been shown that the solution has a positive lower bound $u_\epsilon(h) > \epsilon$. However, the main obstacle remains as the perturbation approaches zero. Since after taking the limit $\epsilon \to 0$ it is only guaranteed that $u(h) > 0$ almost everywhere, we can not prevent $u$ from touching down to zero where $u^{-1}$ tends to infinity. More explicitly, $1/u_\epsilon$ is in $L^1$ space, which is non-reflexive, and there is no weak compactness in $L^1$ space. Thus as we take the limit $\epsilon \to 0$, $\frac{1}{u}$ is in Radon space $\mathcal{M}([0,1])$, rather than $L^1$ space. Hence we can not obtain $\left(\frac{1}{u}\right)_t = (u^3)_{hhhh}$ by directly taking the limit in the regularized problem.

In order to improve the results for (1.6) and explore other interesting dynamics in solid films, we modify the surface energy (1.3) to incorporate a logarithmic factor,

$$G(h) = \int \alpha |h_x| \ln |h_x| + \tfrac{1}{2}|h_x|^3 \, \mathrm{d}x. \tag{1.9}$$

This modification is motivated by kinetic theory and related energy techniques, and the contribution from the new energy enables us to gain weak compactness in the proof of global strong solutions. This modified surface energy is comparable to (1.3) since the logarithmic correction is negligible for small surface gradients.

With this new energy, the one-dimensional evolution equation, restricted to the case when $h_x > 0$ on $0 \le x \le L$, is

$$h_t = -\left[\frac{1}{h_x}\left(\alpha \ln(h_x) + \tfrac{3}{2}h_x^2\right)_{xx}\right]_x, \tag{1.10}$$

This evolution equation is comparable to (1.5) except for the logarithmic term due to the difference between the energy functionals (1.3) and (1.9). The rate of dissipation of the surface energy (1.9) for this model is

$$\frac{\mathrm{d}G}{\mathrm{d}t} = \int_0^L h_t \frac{\delta G}{\delta h} \, \mathrm{d}x = -\int_0^L \frac{1}{h_x}\left(\frac{\delta G}{\delta h}\right)_x^2 \mathrm{d}x \le 0.$$

Applying the change of variables $u(t,h) = h_x(t,x)$ to (1.10) yields the slope equation

$$u_t = -u^2(u^3 + \alpha u)_{hhhh}. \tag{1.11}$$

This one-dimensional fourth-order nonlinear PDE (1.11) is the main focus of this paper. Note that the previously studied model (1.6) corresponds to the case $\alpha = 0$ in (1.11). Moreover, the continuum model (1.11) can also be derived from (1.8) if the same $r^{-2}$ and $\ln|r|$ elastic contributions in [8, 26] with only local step interactions are considered, with the local contributions in (1.8) set to be

$$f(r) = \frac{1}{2r^2} - \alpha \ln |r|. \tag{1.12}$$

This form of $f(r)$ is consistent with the choice of the surface energy in [26] where elastic interactions among steps are incorporated in the derivation of chemical potentials from a discrete BCF model.

For the case $\alpha > 0$, we introduce the scalings

$$\tilde{u} = \alpha^{-1/2} u, \qquad \tilde{h} = \alpha^{-1/2} h, \qquad \tilde{t} = t$$

for the governing slope equation (1.11), which after dropping the tildes leads to

$$u_t = -u^2(u^3 + u)_{hhhh}. \tag{1.13}$$

This rescaling simplifies the discussion of the original problem (1.11).

To model a monotone step train with periodic slope where the maximum and minimum heights in each period do not change in time, we impose periodic boundary conditions on $u$,

$$u(h) = u(h + H), \tag{1.14}$$

corresponding to Dirichlet boundary conditions on $h$ with a fixed height difference $H$ in each period,

$$h_x(x) = h_x(x + L), \qquad h(0) = 0, \qquad h(L) = H. \tag{1.15}$$

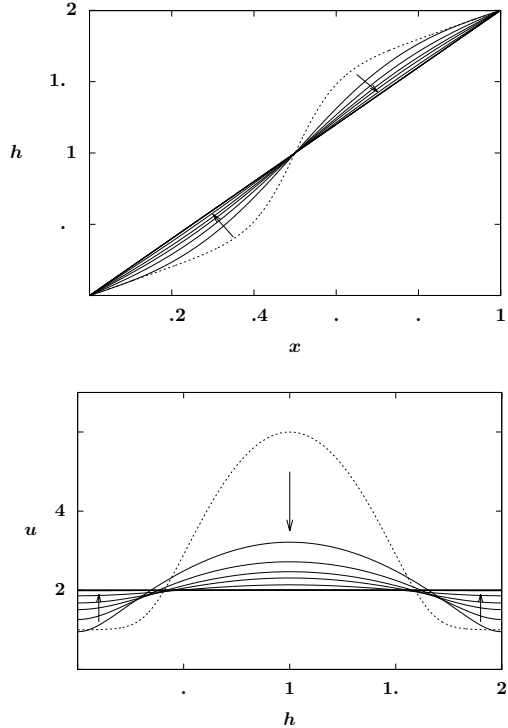

FIGURE 1. (Top) A typical PDE simulation for (1.10) with $\alpha = 1$ on $0 \leq x \leq 1$ and (bottom) the corresponding plot of $u(t, h)$ with boundary conditions (1.15) and $H = 2$, clearly showing the convergence of $h$ to a straight line, with the slope $u$ approaching to a spatially-uniform profile $u = 2$.

Specifically, we investigate the regularity of solutions of (1.13) associated with the periodic boundary condition (1.14) and positive initial data

$$u_0(h) = u(0, h) \geq c_0 > 0, \quad \text{for some constant } c_0 > 0. \tag{1.16}$$

In this work, we obtain a strictly positive lower bound for $u$, which prevents $1/u$ from blowing up, and prove the existence and uniqueness of global strong solutions to (1.13) with periodic boundary conditions. We also study the long time behavior of $u$ by investigating two related Lyapunov functionals for (1.13). In the long time, $u(t, h)$ converges to a spatially-uniform solution which only depends on the initial data $u_0$. Fig. 1 shows typical evolution of the $h$-equation (1.10) solved numerically using finite differences with $\alpha = 1$, $L = 1$ and $H = 2$ starting from the monotone initial data $h_0(x) = 0.5 \tanh(10(x - 0.5)) + 200(x - 0.5) + 1$ on $0 \leq x \leq 1$. It is observed that the monotone profile of $h$ converges to the line $h = 2x$ as $t \to \infty$, with the corresponding slope profile $u$ approaching to a spatially-uniform solution $u^\star = 2$ as $t \to \infty$.

Since the model (1.13) for $0 \leq h \leq H$ is equivalent to the model with $0 \leq h \leq 1$ up to a rescaling, for the rest of this paper we set $H = 1$.

The structure of the rest of the paper is as follows. In section 2 we introduce two useful Lyapunov functions. In section 3 we show the regularity of solutions of PDE (1.13) and the long time behavior of the PDE solution is further investigated in section 4. Numerical verification of these analytical results are then provided in section 5 with a brief discussion of the effect of the logarithmic term on the transient dynamics of the slope equation.

2. **Two Lyapunov functions.** It is important to note that if $u$ is strictly positive, then (1.13) can be written as $(u^{-1})_t = (u^3 + u)_{hhhh}$, which yields the conservation law,
$$\int_0^1 \frac{1}{u}\,dh = \int_0^1 \frac{1}{u_0}\,dh. \tag{2.1}$$
In particular, using the relationship $u(t,h) = h_x(t,x)$ and boundary conditions (1.14) and (1.15), we conclude that
$$\int_0^1 \frac{1}{u}\,dh = \int_0^1 \frac{\partial x}{\partial h}\,dh = x|_{h=1} - x|_{h=0} = L, \tag{2.2}$$
which gives an underlying connection between the surface height equation (1.10) and the slope equation (1.13).

Two Lyapunov energy functions of (1.13) will now be introduced. Specifically, there is a special bi-variational structure embedded in the slope equation (1.13) that is useful in the proof of the regularity of solutions. With the energy density function $f$ in (1.12), the first Lyapunov function is given by
$$F(u) := \int_0^1 \tfrac{1}{2}u^2 + \ln u\,dh. \tag{2.3}$$
Then the variational structure of (1.13) can be written as
$$u_t = -u^2\Big(u^2 \frac{\delta F}{\delta u}\Big)_{hhhh},$$
which leads to the dissipation inequality for $F(u)$,
$$\frac{dF}{dt} = \int_0^1 u_t \frac{\delta F}{\delta u}\,dh = -\int_0^1 (u^3+u)(u^3+u)_{hhhh}\,dh = -\int_0^1 (u^3+u)_{hh}^2\,dh \leq 0. \tag{2.4}$$
In addition, we define the Lyapunov function
$$E(u) := \int_0^1 \big((u^3+u)_{hh}\big)^2\,dh, \tag{2.5}$$
and direct calculation gives that $\delta E/\delta u = (2+6u^2)(u^3+u)_{hhhh}$. Hence the PDE (1.13) actually has a "bi-variational structure" which is given by
$$u_t = -\frac{u^2}{2+6u^2}\frac{\delta E}{\delta u} = -u^2\Big(u^2\frac{\delta F}{\delta u}\Big)_{hhhh}. \tag{2.6}$$
This is the key point in proving the existence and positivity of $u$. From (2.6), we also obtain the dissipation inequality for $E(u)$,
$$\frac{dE}{dt} = \int_0^1 u_t \frac{\delta E}{\delta u}\,dh = -\int_0^1 \frac{u^2}{2+6u^2}\Big(\frac{\delta E}{\delta u}\Big)^2\,dh \leq 0. \tag{2.7}$$
Moreover, from (2.4) and (2.5), the relation between $F(u)$ and $E(u)$ is
$$\frac{dF}{dt} + E = 0. \tag{2.8}$$

This, together with (2.7), gives that

$$TE(T) \leq \int_0^T E(u(t))\,\mathrm{d}t = F(u(0)) - F(u(T)) \leq F(u(0)) - \int_0^1 \ln u_{\min}\,\mathrm{d}h, \quad (2.9)$$

and thus

$$E(u(T)) \leq \frac{C_0}{T} \qquad \text{for any } T > 0,$$

where $C_0$ is a constant depending only on initial data $u_0$, and $u_{\min}$ is the minimal value of $u$ over $[0,T] \times [0,H]$. This formal observation is important for studying the long time behavior of strong solutions to (1.13).

From the energy dissipation (2.7), to obtain an exponential decay for energy $E$, it remains to prove that the dissipation rate of $E$ should be bounded above by $-E$, i.e., $\frac{\mathrm{d}E}{\mathrm{d}t} \leq -cE$, which is usually shown by a logarithmic Sobolev inequality. One example for this is [25]. After proving the global positive lower bound for $u$, we will use Poincare's inequality to obtain this relation and an exponential decay rate for the energy $E$. See the analysis for the long time behavior of strong solutions in section 4.

3. **Global strong solution.** In the following, with standard notations for Sobolev spaces, we denote

$$H^m_{\mathrm{per}}([0,1]) := \{u(h) \in H^m(\mathbb{R}); \ u(h+1) = u(h) \ a.e. \ \in \mathbb{R}\}. \quad (3.1)$$

First we give the definition of a strong solution to PDE (1.13).

**Definition 1.** For any $T > 0$, we define a strong solution $u(t)$ to PDE (1.13) to be a positive function that satisfies

$$u \in C([0,T]; H^2_{\mathrm{per}}([0,1])) \cap L^2([0,T]; H^4_{\mathrm{per}}([0,1])), \quad (3.2)$$

$$u_t \in L^2([0,T]; L^2_{\mathrm{per}}([0,1])) \quad (3.3)$$

with initial data $u_0$ and

$$u_t = -(u^3 + u)_{hhhh} \quad \text{for a.e. } (t,h) \in [0,T] \times [0,1]. \quad (3.4)$$

We now state the main result, the global existence of strong solutions to (1.13) as follows.

**Theorem 1.** *For any $T > 0$, assume initial data $u_0 \in H^2_{per}([0,1])$, $\int_0^1 \frac{1}{u_0}\,\mathrm{d}h = L$ and $u_0 \geq 0$. Then there exists a unique global positive strong solution to PDE (1.13) with initial data $u_0$ and the following two energy-dissipation equalities hold*

$$E(u(T,\cdot)) + \int_0^T \int_0^1 2(u^2 + 3u^4)((u^3+u)_{hhhh})^2\,\mathrm{d}h\,\mathrm{d}t = E(u(0,\cdot)), \quad (3.5)$$

$$F(u(T,\cdot)) + \int_0^T E(u(t,\cdot))\,\mathrm{d}t = F(u(0,\cdot)). \quad (3.6)$$

*Further, the lower bound for $u$ is obtained*

$$u(t,h) \geq \min\left\{\frac{1}{3L^3 E(t)}, \frac{1}{2L} - \frac{E(t)^{\frac{1}{2}}}{3\sqrt{2}}\right\} \geq \min\left\{\frac{1}{3L^3 E_0}, \frac{1}{6L}\right\} > 0, \quad (3.7)$$

*for $(t,h) \in [0,T] \times [0,1]$, where $E_0 = E(u(0))$.*

We also point out the following two important facts:

**Remark 1.** If it is assumed that the initial data $u_0$ is smooth in addition to the conditions in Theorem 1, then with the lower bound in (3.7) we can use standard arguments to obtain higher-order estimates for $u$, and this therefore ensures positive smooth solutions to PDE (1.13).

**Remark 2.** Considering the limit $E(t) \to 0$ as $t \to \infty$ in (3.7) in Theorem 1, we will have an asymptotic lower bound for $u$,

$$u(t) \geq \frac{1}{2L} \qquad \text{for } t \to \infty. \tag{3.8}$$

First we introduce two lemmas which will be used later.

**Lemma 1.** *For any 1-periodic function $v(h)$, we have the following relation*

$$\int_0^1 ((v^3)_{hh})^2 \, dh = 9 \int_0^1 v^4 (v_{hh})^2 \, dh. \tag{3.9}$$

*Proof.* Notice that

$$\begin{aligned}
((v^3)_{hh})^2 &= [(3v^2 v_h)_h]^2 = [6v v_h^2 + 3v^2 v_{hh}]^2 \\
&= 9v^4 v_{hh}^2 + 36v^2 v_h^4 + 36v^3 v_h^2 v_{hh} \\
&= 9v^4 v_{hh}^2 + 36v^2 v_h^4 + 12v^3 (v_h^3)_h \\
&= 9v^4 v_{hh}^2 + 12(v^3 v_h^3)_h.
\end{aligned}$$

Integrating from 0 to 1, we obtain (3.9). □

**Lemma 2.** *For any function $v(h)$ such that $v_{hh} \in L^2([0,1])$, assume that $v$ achieves its minimal value at $h^\star$, i.e. $v_{\min} = v(h^\star)$. Then we have*

$$v(h) - v_{\min} \leq \tfrac{2}{3} \|v_{hh}\|_{L^2([0,1])} |h - h^\star|^{\frac{3}{2}}. \tag{3.10}$$

*Proof.* Since $v_{hh} \in L^2([0,1])$, $v_h$ is continuous. Hence by $v_{\min} = v(h^\star)$, we have $v_h(h^\star) = 0$ and

$$|v_h(h)| = \left| \int_{h^\star}^1 v_{hh}(s) \, ds \right| \leq |h - h^\star|^{\frac{1}{2}} \|v_{hh}\|_{L^2([0,1])}. \tag{3.11}$$

Hence we have

$$\begin{aligned}
|v(h) - v_{\min}| &\leq \int_{h^\star}^1 |s - h^\star|^{\frac{1}{2}} \|v_{hh}\|_{L^2([0,1])} \, ds \\
&\leq \tfrac{2}{3} |h - h^\star|^{\frac{3}{2}} \|v_{hh}\|_{L^2([0,1])}.
\end{aligned}$$

□

*Proof of Theorem 1.* Recall that $L = \int_0^1 \frac{1}{u_0} \, dh$ from (2.2), and denote $u_{\min}$ as the minimal value of $u$ in $[0,T] \times [0,1]$. For any $t \in [0,T]$, we also denote $u_m(t)$ as the minimal value of $u(t)$ for $h \in [0,1]$. Assume $u$ achieves its minimal value at $t^\star$, $h^\star$, i.e. $u_{\min} = u(t^\star, h^\star)$. Notice

$$E_0 = E(u_0) = \int_0^1 [(u_0^3 + u_0)_{hh}]^2 \, dh \leq C(\|u_0\|_{H^2([0,1])}),$$

due to Lemma 1. In Step 1, we first test some estimates under the a priori assumption

$$u(t,h) \geq u_{\min} \geq \frac{1}{2} \min \left\{ \frac{1}{3L^3 E_0}, \frac{1}{6L} \right\} > 0 \qquad \text{for any } t \in [0,T], \ h \in [0,1]. \tag{3.12}$$

In Step 2, we will show that this a priori assumption leads to a better lower bound. Then we can take the limit using those a priori estimates in Step 1. In Step 3, we prove the energy-dissipation equalities (3.5), (3.6) and thus obtain the existence result to (1.13). In Step 4, we prove that the solution obtained above is unique.

**Step 1**: A priori estimates.

First, we obtain the higher order estimate. Multiplying (1.13) by $(u^3 + u)_{hhhh}$ and integrating by parts leads to

$$\int_0^1 (u^3)_{hh} u_{hht} + \frac{1}{2}((u_{hh})^2)_t \, dh$$
$$= \int_0^1 (u^3 + u)_{hhhh} u_t \, dh = -\int_0^1 u^2((u^3 + u)_{hhhh})^2 \, dh \leq 0. \quad (3.13)$$

Then multiplying (1.13) by $3u^2(u^3 + u)_{hhhh}$ and integrating by parts yields

$$\int_0^1 (u^3)_{hht} u_{hh} + \tfrac{1}{2}((u_{hh}^3)^2)_t \, dh$$
$$= \int_0^1 (u^3 + u)_{hhhh}(u^3)_t \, dh = -\int_0^1 3u^4((u^3 + u)_{hhhh})^2 \, dh \leq 0. \quad (3.14)$$

Combining (3.13) and (3.14), we have

$$\frac{dE}{dt} = \frac{d}{dt}\int_0^1 ((u^3)_{hh} + u_{hh})^2 \, dh = -\int_0^1 2(u^2 + 3u^4)((u^3 + u)_{hhhh})^2 \, dh \leq 0. \quad (3.15)$$

Then for any $T > 0$,

$$\int_0^1 (u_{hh} + (u^3)_{hh})^2 \, dh - E_0 = -\int_0^T \int_0^1 2(u^2 + 3u^4)((u^3 + u)_{hhhh})^2 \, dh \, dt \leq 0, \quad (3.16)$$

which gives

$$\|(u + u^3)_{hh}\|_{L^2([0,1])} = E^{\frac{1}{2}} \leq E_0^{\frac{1}{2}}. \quad (3.17)$$

From (3.16) and (3.12), we also have

$$\|(u^3 + u)_{hhhh}\|_{L^2([0,T];L^2([0,1]))} \leq C(\|u_0\|_{H^2([0,1])}, L). \quad (3.18)$$

Second, we use the a priori assumption (3.12) to get the lower order estimate. Denote $\langle u \rangle := \int_0^1 u \, dh$ and $\bar{u} := u - \langle u \rangle$. From Poincare's inequality and (3.17), we have

$$\|u_h\|_{L^2([0,1])}^2 \leq \|u_h\|_{L^2([0,1])}^2 + 9\|u^2 u_h\|_{L^2([0,1])}^2 + 6\|u_h u\|_{L^2([0,1])}^2$$
$$= \|(u^3 + u)_h\|_{L^2([0,1])}^2 \leq E \leq E_0,$$

for any $t \in [0, T]$. Then by Poincare's inequality we obtain

$$\int_0^1 |u - \langle u \rangle|^2 \, dh = \int_0^1 |\bar{u}|^2 \, dh \leq c\|u_h\|_{L^2([0,1])}^2 \leq cE \leq cE_0, \quad (3.19)$$

and it remains to estimate $\langle u \rangle$. Using Poincare's inequality and Sobolev embedding again, we have

$$\|\bar{u}\|_{L^\infty} = \|u - \langle u \rangle\|_{L^\infty([0,1])} \leq \|u - \langle u \rangle\|_{H^1([0,1])}$$
$$\leq c\|(u - \langle u \rangle)_h\|_{L^2([0,1])} \leq cE^{\frac{1}{2}} \leq cE_0^{\frac{1}{2}}. \quad (3.20)$$

On the other hand, since $u$ has the lower bound assumption (3.12), we have (2.1) and (2.2). Hence

$$L = \int_0^1 \frac{1}{u} \, dh = \int_0^1 \frac{1}{\langle u \rangle + \bar{u}} \, dh \leq \int_0^1 \frac{1}{\langle u \rangle + \bar{u}_{\min}} \, dh, \tag{3.21}$$

which implies

$$\langle u \rangle \leq \frac{1}{L} - \bar{u}_{\min} \leq \frac{1}{L} + \|\bar{u}\|_{L^\infty}. \tag{3.22}$$

This, together with (3.20), yields

$$\langle u \rangle \leq \frac{1}{L} + cE^{1/2} \leq \frac{1}{L} + cE_0^{1/2}.$$

Thus for any $T > 0$ we have the lower order estimate

$$\|u\|_{L^\infty([0,T];L^2[0,1])} \leq C(\|u_0\|_{H^2([0,1])}, L), \tag{3.23}$$

which also leads to

$$u_m(t) := \min_{h \in [0,1]} u(t) \leq \|u\|_{L^\infty([0,T];L^2[0,1])} \leq C(\|u_0\|_{H^2([0,1])}, L) \text{ for any } t \in [0, T]. \tag{3.24}$$

Third, notice that the function $u^3 + u$ is increasing. From Lemma 2, we have

$$[u^3(t) + u(t)] - [u_m^3(t) + u_m(t)] \leq \tfrac{2}{3}\|(u^3 + u)_{hh}(t, \cdot)\|_{L^2([0,1])}|h - h^\star|^{\tfrac{3}{2}} \text{ for any } t \in [0, T]. \tag{3.25}$$

This, together with (3.17) and (3.24), shows that

$$\|u^3 + u\|_{L^\infty([0,T] \times [0,1])} \leq C(\|u_0\|_{H^2([0,1])}, L). \tag{3.26}$$

Hence we have

$$\|u\|_{L^\infty([0,T] \times [0,1])} \leq C(\|u_0\|_{H^2([0,1])}, L), \tag{3.27}$$

and

$$\|u^3 + u\|_{L^\infty([0,T];H^2[0,1])} \leq C(\|u_0\|_{H^2([0,1])}, L). \tag{3.28}$$

Next, we define $v = \varphi(u) := u^3 + u$. From a priori assumption (3.12) and the estimate (3.27), we know that $u = \varphi^{-1}(v)$ is smooth on $[v_{\min}, v_{\max}]$. Hence from Sobolev embedding we have

$$\|u(t, \cdot)\|_{H^2([0,1])} = \|\varphi^{-1}(v(t, \cdot))\|_{H^2([0,1])} \leq C(\|v(t, \cdot)\|_{L^\infty([0,1])})\|v(t, \cdot)\|_{H^2([0,1])}$$

for any $t \in [0, T]$. This, together with (3.27) and (3.28), implies the estimate

$$\|u\|_{L^\infty([0,T];H^2([0,1]))} \leq C(\|u_0\|_{H^2([0,1])}, L). \tag{3.29}$$

and $u \in L^\infty([0, T]; H^2([0, 1]))$.

Similarly, we have

$$\|u(t, \cdot)\|_{H^4([0,1])} = \|\varphi^{-1}(v(t, \cdot))\|_{H^4([0,1])} \leq C(\|v(t, \cdot)\|_{L^\infty([0,1])})\|v(t, \cdot)\|_{H^4([0,1])}$$

for any $t \in [0, T]$. This, together with (3.27) and (3.18), implies that $u \in L^2([0, T]; H^4([0, 1]))$.

Finally, we can also get the estimate for $u_t$ from (3.15),

$$\int_0^T \int_0^1 u_t^2 \, dh \, dt = \int_0^T \int_0^1 u^4((u^3 + u)_{hhhh})^2 \, dh \, dt \leq E_0.$$

Hence we have

$$\|u_t\|_{L^2([0,T];L^2([0,1]))} \leq E_0^{\tfrac{1}{2}}. \tag{3.30}$$

Moreover, since $u \in L^2([0,T]; H^4([0,1]))$ and $u_t \in L^2([0,T]; L^2([0,1]))$ by [6, Theorem 4, p. 288], we have $u \in C([0,T]; H^1([0,1]))$ and for any $0 \leq s \leq \tau \leq T$

$$\int_s^\tau \frac{d}{dt} \|u_{hh}\|_{L^2([0,1])}^2 \, dt$$
$$= \int_s^\tau \int_0^1 2u_t \partial_h^4 u \, dh \, dt \leq \|u_t\|_{L^2([s,\tau]; L^2([0,1]))} + \|\partial_h^4 u\|_{L^2([s,\tau]; L^2([0,1]))},$$

which implies $u \in C([0,T]; H^2([0,1]))$.

**Step 2:** Verify the a priori assumption.

From (3.25), we have

$$[u(t) - u_m(t)] \leq \frac{2}{3} E(t)^{\frac{1}{2}} |h - h^\star|^{\frac{3}{2}} \qquad \text{for any } t \in [0,T]$$

and for any $0 < \beta \leq \frac{1}{2}$, using the relationship in (2.2) we obtain

$$\frac{\beta}{u_m(t) + \frac{2E(t)^{\frac{1}{2}}}{3} \beta^{\frac{3}{2}}} = \int_{h^\star}^{h^\star + \beta} \frac{1}{u_m(t) + \frac{2E(t)^{\frac{1}{2}}}{3} \beta^{\frac{3}{2}}} \, dh$$
$$\leq \int_0^1 \frac{1}{u_m(t) + \frac{2E(t)^{\frac{1}{2}}}{3} |h - h^\star|^{\frac{3}{2}}} \, dh$$
$$\leq \int_0^1 \frac{1}{u} \, dh = L.$$

Hence we have

$$u_m(t) \geq \frac{\beta}{L} - \frac{2E^{\frac{1}{2}}}{3} \beta^{\frac{3}{2}},$$

and for $0 < \beta \leq \frac{1}{2}$ the right hand side attains its maximum value at

$$\beta = \begin{cases} \frac{1}{2} & \text{for } \frac{1}{E(t)L^2} > \frac{1}{2}, \\ \frac{1}{E(t)L^2} & \text{for } \frac{1}{E(t)L^2} \leq \frac{1}{2}. \end{cases}$$

Then direct calculation leads to a lower bound $u_m$ that depends on $E(t)$

$$u_m(t) \geq \mathcal{J}(E(t)), \tag{3.31}$$

where

$$\mathcal{J}(E(t)) = \begin{cases} \dfrac{1}{3L^3 E(t)} & \text{for } E(t) \geq \dfrac{2}{L^2}, \\ \dfrac{1}{2L} - \dfrac{E(t)^{\frac{1}{2}}}{3\sqrt{2}} & \text{for } E(t) < \dfrac{2}{L^2}. \end{cases} \tag{3.32}$$

Therefore we have

$$u_m(t) \geq \min\left\{\frac{1}{3L^3 E(t)}, \frac{1}{2L} - \frac{E(t)^{\frac{1}{2}}}{3\sqrt{2}}\right\} \qquad \text{for any } t \in [0,T], \tag{3.33}$$

and taking the minimal value with respect to $t \in [0,T]$, we have

$$u_{\min} \geq \min\left\{\frac{1}{3L^3 E_0}, \frac{1}{6L}\right\} > \frac{1}{2} \min\left\{\frac{1}{3L^3 E_0}, \frac{1}{6L}\right\}, \tag{3.34}$$

which verifies the a priori assumption (3.12) and shows that $u_{\min}$ has a positive lower bound,

$$u_{\min} \geq \min\left\{\frac{1}{3L^3 E_0}, \frac{1}{6L}\right\} > 0.$$

Thus we obtain (3.7).

**Step 3:** Proof of energy-dissipation equalities.

After obtaining the a priori estimates in Step 1, using standard compactness argument we obtain that for any $T > 0$ and any $\phi_i \in C^\infty([0,T] \times [0,1])$, $u$ satisfies

$$\int_0^T \int_0^1 \phi_i u_t \, dh \, dt + \int_0^T \int_0^1 \phi_i u^2 (u^3 + u)_{hhhh} \, dh \, dt = 0. \tag{3.35}$$

This, together with the fundamental lemma of the calculus of variations, yields (3.4). Since we have a positive lower bound (3.7) for $u$, we can take $\phi_i$ such that $\phi_i \to u + \frac{1}{u}$ in $L^2([0,T]; L^2([0,1]))$ as $i \to \infty$. Hence from (3.2) and (3.3), we take a limit in (3.35) to obtain

$$\int_0^T \int_0^1 \left(\tfrac{1}{2}u^2 + \ln u\right)_t dh \, dt = -\int_0^T \int_0^1 (u^3 + u)(u^3 + u)_{hhhh} \, dh \, dt$$
$$= -\int_0^T \int_0^1 ((u^3 + u)_{hh})^2 \, dh \, dt = -\int_0^T E \, dt, \tag{3.36}$$

which implies (3.6).

Similarly, since (3.2) implies $(1 + 3u^2)(u^3 + u)_{hhhh} \in L^2([0,T]; L^2([0,1]))$, we can take $\phi_i$ in (3.35) such that $\phi_i \to 2(1 + 3u^2)(u^3 + u)_{hhhh}$ in $L^2([0,T]; L^2([0,1]))$ as $i \to \infty$. Hence from (3.2) and (3.3) again, we take a limit in (3.35) to obtain

$$\int_0^T \frac{dE}{dt} \, dt = \int_0^T \int_0^1 \left((u^3 + u)_{hh}^2\right)_t dh \, dt$$
$$= -\int_0^T \int_0^1 2(3u^4 + u^2)((u^3 + u)_{hhhh})^2 \, dh \, dt, \tag{3.37}$$

which implies (3.5).

**Step 4:** Uniqueness of the solution to (1.13).

Assume that $u, v$ are two solutions of (1.13). Then we have

$$(u - v)_t = -u^2(u^3 + u)_{hhhh} + v^2(v^3 + v)_{hhhh}, \tag{3.38}$$

$$(u^3 - v^3)_t = -3u^4(u^3 + u)_{hhhh} + 3v^4(v^3 + v)_{hhhh}. \tag{3.39}$$

Combining (3.39) and (3.38), we have

$$[u^3 + u - v^3 - v]_t = (3v^2 + 1)v^2(v^3 + v)_{hhhh} - (3u^2 + 1)u^2(v^3 + v)_{hhhh} \tag{3.40}$$
$$+ (3u^2 + 1)u^2(v^3 + v)_{hhhh} - (3u^2 + 1)u^2(u^3 + u)_{hhhh}.$$

Recall that for any $p \geq 0$, $u^p$ is increasing with respect to $u$, so from (3.7) there exist constants $m, M > 0$ whose values depend on $\|u_0\|_{H^2([0,1])}$, $p$ and $L$ and satisfy

$$m \leq u^p \leq M, \qquad m \leq v^p \leq M. \tag{3.41}$$

First, multiplying (3.40) by $(u^3 + u - v^3 - v)_{hhhh}$ and integrating by parts, from Young's inequality and (3.18) we obtain

$$\int_0^1 \left[(3v^4 - 3u^4 + v^2 - u^2)(v^3 + v)_{hhhh}\right](u^3 + u - v^3 - v)_{hhhh}\,\mathrm{d}h$$
$$\leq c(m,M)\|u-v\|_{L^\infty} + \frac{c(m)}{4}\int_0^1 [(u^3 + u - v^3 - v)_{hhhh}]^2 \,\mathrm{d}h,$$

and

$$\frac{\mathrm{d}}{\mathrm{d}t}\int_0^1 \tfrac{1}{2}[(u^3 + u - v^3 - v)_{hh}]^2 \,\mathrm{d}h$$
$$= \int_0^1 \left[(3v^4 - 3u^4 + v^2 - u^2)(v^3 + v)_{hhhh}\right](u^3 + u - v^3 - v)_{hhhh}\,\mathrm{d}h$$
$$- \int_0^1 (3u^2 + 1)u^2[(u^3 + u - v^3 - v)_{hhhh}]^2 \,\mathrm{d}h$$
$$\leq c(m,M)\|u-v\|_{L^\infty} + \frac{c(m)}{4}\int_0^1 [(u^3 + u - v^3 - v)_{hhhh}]^2 \,\mathrm{d}h$$
$$- c(m)\int_0^1 [(u^3 + u - v^3 - v)_{hhhh}]^2 \,\mathrm{d}h$$
$$\leq c(m,M)\|u-v\|_{L^\infty} - \frac{3c(m)}{4}\int_0^1 [(u^3 + u - v^3 - v)_{hhhh}]^2 \,\mathrm{d}h, \tag{3.42}$$

where (3.41) is used, and $c(m,M)$ is a constant that only depends on $m$, $M$.

Second, we multiply (3.40) by $u^3 + u - v^3 - v$ and use integration by parts again. From (3.18), (3.41) and Young's inequality, we have

$$\frac{\mathrm{d}}{\mathrm{d}t}\int_0^1 \tfrac{1}{2}[u^3 + u - v^3 - v]^2 \,\mathrm{d}h$$
$$= \int_0^1 \left[(3v^4 - 3u^4 + v^2 - u^2)(v^3 + v)_{hhhh}\right](u^3 + u - v^3 - v)\,\mathrm{d}h \tag{3.43}$$
$$- \int_0^1 (3u^2 + 1)u^2 (u^3 + u - v^3 - v)_{hhhh}(u^3 + u - v^3 - v)\,\mathrm{d}h$$
$$\leq c(m,M)\|u-v\|_{L^2([0,1])} + \frac{c(m)}{4}\int_0^1 [(u^3 + u - v^3 - v)_{hhhh}]^2 \,\mathrm{d}h.$$

Finally, (3.42) and (3.43) show that

$$\frac{\mathrm{d}}{\mathrm{d}t}\|u^3 + u - v^3 - v\|_{H^2([0,1])}^2 \leq c(m,M)\|v-u\|_{L^\infty([0,1])}^2. \tag{3.44}$$

Again from (3.7), we have

$$|u - v| \leq |u - v|(u^2 + v^2 + uv + 1) = |u^3 - v^3 + u - v|.$$

Therefore

$$\frac{\mathrm{d}}{\mathrm{d}t}\|u^3 + u - v^3 - v\|_{H^2([0,1])}^2 \leq c(m,M)\|v-u\|_{L^\infty([0,1])}^2$$
$$\leq c(m,M)\|u^3 + u - v^3 - v\|_{L^\infty([0,1])}^2 \tag{3.45}$$
$$\leq c(m,M)\|u^3 + u - v^3 - v\|_{H^2([0,1])}^2.$$

Hence if $u(0) = v(0)$, Grönwall's inequality implies that $u = v$. This completes the proof of Theorem 1. □

4. **Long time behavior of strong solution.** After establishing the global-in-time strong solution, we want to study how the solution will behave for long times. For the PDE (1.13) with periodic boundary conditions, the solution converges to a constant.

**Theorem 2.** *Under the same assumptions of Theorem 1, given the strong solution $u$ obtained in Theorem 1, we have an exponentially decay estimate for energy $E$, i.e.*

$$E(u(t)) \leq E_0 e^{-\beta t}, \ as \ t \to +\infty, \tag{4.1}$$

*where $\beta$ is a positive constant depending on $u_0$. Moreover, there exists a constant $u^\star = 1/L$ such that, as time $t \to \infty$, $u$ converges to $u^\star$ in the sense*

$$\|u(t, \cdot) - u^\star\|_{L^\infty([0,1])} \to 0, \ as \ t \to \infty. \tag{4.2}$$

*Proof.*
**Step 1:** Exponential decay of the free energy $E$ as $t \to \infty$.

By (3.15) we know that $E$ is decreasing with respect to $t$. Combining (3.15) and the lower bound estimate (3.7), we have

$$\frac{\mathrm{d}}{\mathrm{d}t} E + (2u_{\min}^2 + 6u_{\min}^4)\|(u^3 + u)_{hhhh}\|_{L^2([0,1])}^2 \leq 0. \tag{4.3}$$

Note that the sharp Poincare's inequality leads to

$$\|(u^3 + u)_{hh}\|_{L^2([0,1])}^2 \leq \|(u^3 + u)_{hhhh}\|_{L^2([0,1])}^2.$$

Thus from (4.3) we obtain

$$\frac{\mathrm{d}}{\mathrm{d}t} E \leq -(2u_{\min}^2 + 6u_{\min}^4)\|(u^3 + u)_{hhhh}\|_{L^2([0,1])}^2$$
$$\leq -(2u_{\min}^2 + 6u_{\min}^4)\|(u^3 + u)_{hh}\|_{L^2([0,1])}^2$$
$$= -(2u_{\min}^2 + 6u_{\min}^4)E = -\beta E$$

with the constant $\beta = 2u_{\min}^2 + 6u_{\min}^4 > 0$. Thus $E(u(t))$ converges to 0 as $t \to \infty$ with an exponential decay rate

$$E(u(t)) \leq E_0 e^{-\beta t},$$

which shows that (4.1) holds.

**Step 2:** Showing convergence to the stationary solution $u^\star = \frac{1}{L}$.

Since we have a positive lower bound for $u$ (3.7), the relations (2.1) and (2.2) hold. Recall the notations in the proof of Theorem 1, $\langle u \rangle = \int_0^1 u \, \mathrm{d}h$ and $\bar{u} = u - \langle u \rangle$. From Hölder's inequality,

$$1 = \left(\int_0^1 \frac{1}{\sqrt{u}} \sqrt{u} \, \mathrm{d}h\right)^2 \leq \int_0^1 u \, \mathrm{d}h \int_0^1 \frac{1}{u} \, \mathrm{d}h,$$

which shows that

$$\frac{1}{L} \leq \frac{1}{\int_0^1 \frac{1}{u} \, \mathrm{d}h} \leq \langle u \rangle.$$

This, together with (3.22), gives

$$\|\langle u \rangle - \frac{1}{L}\|_{L^\infty([0,1])} \leq \|\bar{u}\|_{L^\infty([0,1])}. \tag{4.4}$$

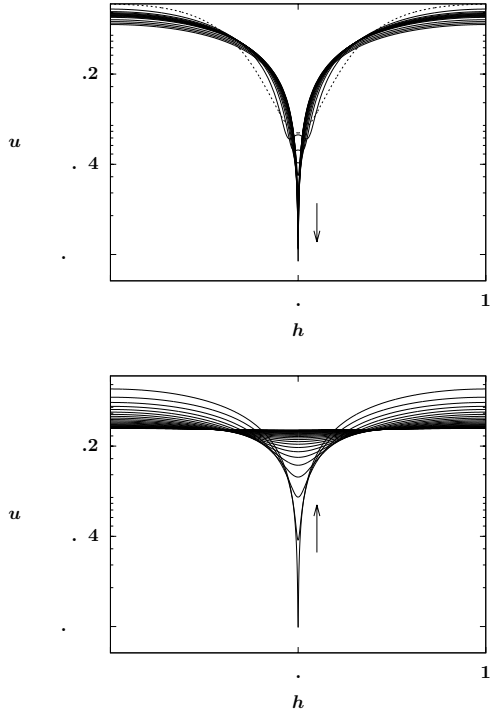

FIGURE 2. A numerical simulation of PDE (1.6) plotted in semi-log coordinates starting from the initial condition (5.1) (plotted with the dashed line): (top) early stage near-rupture is approached as the global minimum decreases from 0.07 to 0.007 for $0 < t < 0.0032$; (bottom) later stage behavior for $t > 0.0032$ as the solution approaches a constant $u^\star = 0.27$.

Thus from (3.20), we obtain

$$\|u - \frac{1}{L}\|_{L^\infty([0,1])} \leq \|u - \langle u \rangle\|_{L^\infty([0,1])} + \|\langle u \rangle - \frac{1}{L}\|_{L^\infty([0,1])} \tag{4.5}$$

$$\leq 2\|\bar{u}\|_{L^\infty([0,1])} \leq cE^{\frac{1}{2}}, \tag{4.6}$$

which, combining with the decay of $E$ in step 1, completes the proof. □

5. **Numerical study.** To simulate the surface growth dynamics and explore beyond the analytical results presented in section 3 and section 4, we numerically investigate both the height equation (1.10) and the slope equations (1.6) and (1.13). While the numerical results presented in Fig. 1 are obtained directly from solving the height equation (1.10) with Dirichlet boundary conditions (1.15) and specified values for $L$ and $H$, similar calculations can be carried out on the slope equations (1.6) and (1.13) with periodic boundary conditions (1.14). Specifically, we are interested in different transient behaviors of the solutions to (1.6) and (1.13) and their long time behavior.

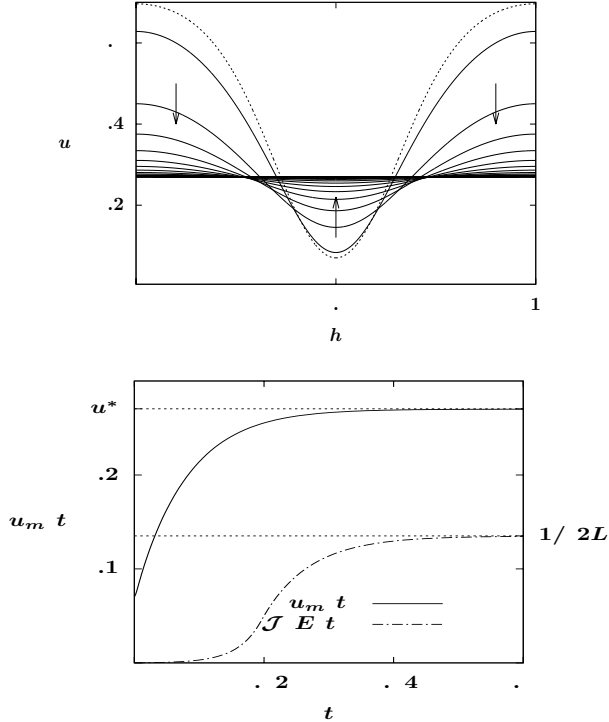

FIGURE 3. (Top) A numerical simulation of PDE (1.13) starting from identical initial conditions used in Fig. 2 showing convergence to a spatially-uniform solution $u = u^\star$ as $t \to \infty$. (Bottom) A plot showing that $u_m(t) = \min_h u(t, h)$ is bounded below by $\mathcal{J}(E(t))$ given by (3.32) which is in line with the conclusion of Theorem 1, and that the asymptotic lower bound $\mathcal{J}(E(t)) \to 1/(2L)$ for $t \to \infty$ as in (3.8).

Typical numerical simulations of (1.6) and periodic boundary conditions are plotted in Figure 2, where the spatial variation grows starting from the initial data

$$u_0(h) = 0.7 - 0.63 \exp(-5(2h-1)^2), \tag{5.1}$$

and leads to a near-rupture transient behavior in the early stage before the solution converges to a constant solution $u = u^\star$ for $t \to \infty$. While for equation (1.13) (or equation (1.11) with $\alpha = 1$), starting from the same initial data, monotone behavior of the PDE solution $u$ converging to $u = u^\star$ is presented in Fig. 3. With the initial condition (5.1), we have $L = \int_0^1 \frac{1}{u_0} \, dh = 3.7$, and the final state $u^\star = 0.27$ can be obtained from the formula $u^\star = 1/L$. The energy-dependent lower bound estimate (3.31) is also plotted in a dashed curve in Fig. 3 (right) in comparison to the numerical results. These numerical studies are conducted using centered finite differences in a Keller box scheme, where the fourth-order PDE (1.13) is decomposed into a system of first-order differential equations

$$u_t + u^2 Q_h = 0, \qquad Q = P_h, \qquad P = K_h, \qquad K = (u^3 + u)_h. \tag{5.2}$$

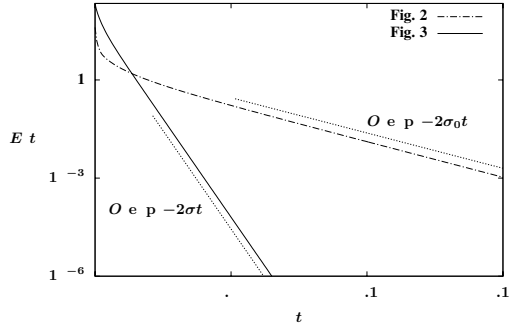

FIGURE 4. Plots of corresponding energy $E$ in (2.5) and (5.7) for PDE simulations in Fig. 2 and Fig. 3. The energy $E(t)$ decays exponentially to zero following (5.6) with $k = 2\pi$.

We now revisit the long-time behavior of the solutions of (1.13). It has been shown in the previous section that as $t \to \infty$, the PDE solution $u(t, h)$ approaches to a spatially-uniform solution $u = u^\star$. To study the long time behavior of the solution, using $u(t, h) = u^\star + \epsilon v(t, h)$ we obtain the linearized equation for (1.13) as $\epsilon \to 0$,

$$v_t = -r(u^\star) v_{hhhh}, \tag{5.3}$$

where the function $r$ is defined as $r(u) = u^2(3u^2 + 1)$. Note that the positivity of solutions to the fourth-order linear PDE (5.3) is not guaranteed. Furthermore, we perturb $\bar{u}$ by individual Fourier mode disturbances

$$u(t, h) = u^\star + \delta e^{ikh} e^{\sigma t} + O(\delta^2), \tag{5.4}$$

where $k$ is the wavenumber and $\sigma$ represents the rate of the PDE solution converging to $u^\star$. Substituting (5.4) into (1.13) and linearizing about $u = u^\star$ leads to the dispersion relation

$$\sigma = -r(u^\star) k^4, \tag{5.5}$$

which indicates that the steady state solution $u \equiv u^\star$ is stable with respect to any Fourier mode perturbations. As the spatially-uniform solution is approached for long times, the energy $E(t)$ decays exponentially in the form of

$$E(t) = C \exp(2\sigma t) = C \exp(-2r(u^\star) k^4 t), \text{ for } t \to \infty. \tag{5.6}$$

where $C$ depends on the initial conditions and other system parameters. We note that the energy decay rate $2r(u^\star) k^4$ obtained from the linear analysis is inline with the estimated bound of energy decay rate $\beta = 2r(u_{\min})$ obtained in Theorem 2. Similar linear analysis for PDE (1.6) leads to the energy decay rate (5.6) where the energy is given by

$$E(u) = \int_0^1 ((u^3)_{hh})^2 \, dh \tag{5.7}$$

and the dispersion relation $\sigma_0 = -3u^4 k^4$. It is shown in Fig. 4 that the exponentially decay rates of the energy $E(t)$ for the PDE simulations in Fig. 2 and Fig. 3 both agree with the above linearization result as $t \to \infty$.

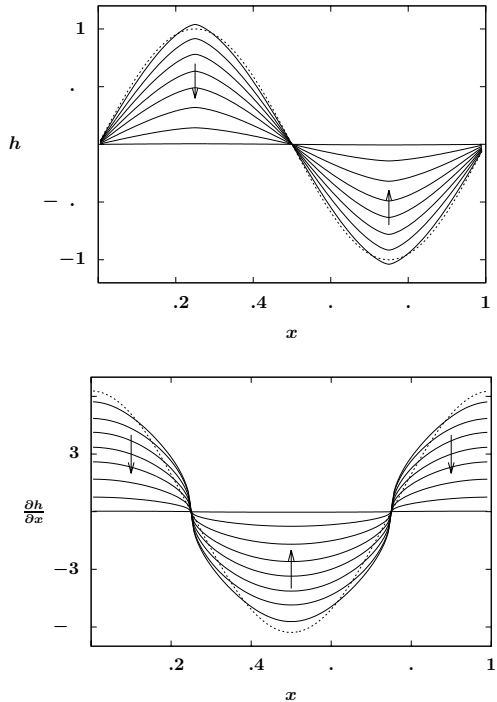

FIGURE 5. Evolution of the surface height $h(t,x)$ and slope $h_x(t,x)$ following equation (6.2) with $\alpha = 0$ starting from initial condition $h_0(x) = \sin(2\pi x)$ on $0 \leq x \leq 1$, showing convergence to spatially-uniform solution $h \equiv 0$ as $t \to \infty$.

6. **Discussion.** The main contribution of this paper is showing the global strong solution for the continuum slope equation (1.13) using the interesting bi-variational structure that the equation possesses. A formal derivation of the slope equation (1.13) from the discrete model with an additional logarithmic energy term is included, as well as its connection with the surface height equation (6.1) under the strict monotonicity assumption. Typical numerical simulations of both the height equation (1.10) and the slope equation (1.13) are also presented in support of the analytical regularity estimate. In particular, we investigate the effect of an additional logarithmic contribution in the slope equation and show that the transient near-rupture behavior occurs in (1.6). As $t \to \infty$, the solution to (1.13) converges exponentially to a constant solution $u^\star$ which is selected by the initial data.

In addition to the evolution of monotone trains of steps with boundary conditions (1.14) and (1.15) described above, we are also inspired by Li and Liu's work on thin film epitaxy [15, 16] to further investigate the surface height evolution $h$ in the periodic setting, which is described by the height equation

$$h_t = -\Big[\frac{1}{|h_x|}\Big(\alpha \frac{h_x}{|h_x|} \ln |h_x| + \tfrac{3}{2}|h_x|h_x\Big)_{xx}\Big]_x \qquad (6.1)$$

with associated periodic boundary conditions $h(t,x) = h(t, x+L)$. Note that in this regime the solution $h$ is allowed to have both positive and negative slopes which is not covered by the argument provided in section 3. While the regularity of solutions of (6.1) in this regime is beyond the scope of this article, we shall present some preliminary numerical results and show the potential connection between (6.1) and thin film epitaxy models with slope selections [15].

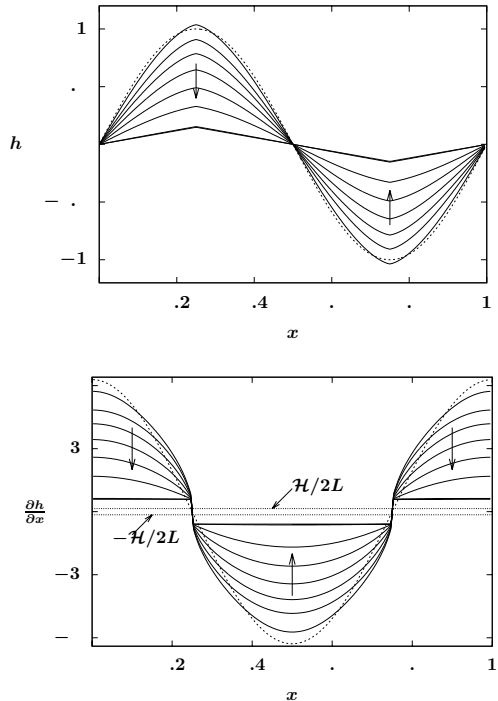

FIGURE 6. Evolution of the surface height $h(t,x)$ and slope $h_x(t,x)$ for equation (6.2) with $\alpha = 1$ starting from identical initial data used in Fig. 5, showing convergence to a piece-wise constant profile in $h$ and jump in $h_x$.

Since the logarithmic term in the chemical potential and the mobility function both become singular as $|h_x| \to 0$, a direct finite-difference discretization of the equation leads to numerical difficulties. Therefore, we introduce regularization for both the logarithmic term and the mobility function in the model with small coefficients $\epsilon$ and $\delta$, and numerically study the following regularized equation,

$$h_t = -\left[\frac{1}{\sqrt{h_x^2 + \epsilon^2}}\left(\alpha \frac{h_x \ln |h_x|}{\sqrt{h_x^2 + \delta^2}} + \tfrac{3}{2}|h_x|h_x\right)_{xx}\right]_x. \tag{6.2}$$

With $\epsilon = \delta = 10^{-5}$, numerical simulations for (6.2) with both $\alpha = 0$ and $\alpha = 1$ starting from identical initial data $h_0(x) = \sin(2\pi x)$ are presented in Fig. 5 and Fig. 6. It is observed that for the case $\alpha = 0$, the height profile flattens out with a shock formed in the slope profile. This is consistent with the numerical results for the two-dimensional height profile without line tension in [24]. However, for $\alpha = 1$, due to the contribution of the logarithmic term in (6.1), a shock is formed in $h$ at the extrema of the height profile, and the solution $h$ approaches to a periodic piece-wise linear function with positive and negative slopes. Although the proof provided in Section 3 and Section 4 cannot be applied to this periodic setting directly, we plot the asymptotic bounds $\pm\mathcal{H}/2L$ for the slope function from (3.8) to compare it with the slope profile obtained from the simulation in Fig. 6, where $\mathcal{H} = 0.2936$ and $L = 1$, and the step height $\mathcal{H}$ is obtained from the difference between the two extrema in the final stage in Fig. 6 (left). It is suggestive that the asymptotic lower bound $\mathcal{H}/2L < |h_x|$ away from the jump in $h_x$ and is still a good estimate of the slope profile.

A number of questions remain to be answered. In the periodic surface height scenario, the interesting dynamics from numerical observation of the regularized problem (6.2) with the non-zero logarithmic contribution need further investigation. While we use classical finite difference methods to numerically simulate the dynamics of the regularized problem (6.2), more sophisticated numerical schemes can be developed to capture the dynamics of the unregularized problem (6.1). For instance, some numerical studies [4, 23, 24] of the well-known equation (1.4) use a nonlinear Galerkin scheme and finite-element schemes have been applied to the corresponding DL regime [13]. In addition, we are also interested to study whether the logarithmic term in the surface energy contribute to the singularity at the edge of the facets and steps in the modeling of crystal surface evolution. Moreover, without the monotone assumption for initial data, the analysis of the general $h$-equations (1.4) and (6.1) is very challenging and still open. Due to the form of the mobility function $|\nabla h|^{-1}$, new techniques are needed to rigorously describe the singularity that occurs at $|\nabla h| = 0$.

Inspired by Xiang's work in [26] which incorporates nonlocal interactions in the DL regime, we are also interested in the nonlocal contributions in the ADL regime both analytically and numerically. Under the monotone assumption, the $h$-equation can be written as

$$h_t = -\left[\frac{1}{h_x}\left(\mathscr{H}(h_x) + \left(\alpha \ln h_x + \frac{3}{2}h_x^2\right)_x\right)_x\right]_x,$$

where $\mathscr{H}$ represents the Hilbert transform. While in this work only the one-dimensional models are investigated, the equations in higher dimensions are also interesting future subjects. Furthermore, since the general PDE (1.1) can be obtained from linearizing the exponential in the Gibbs-Thomson relation by taking the approximation $e^{\frac{\delta G}{\delta h}} \approx 1 + \frac{\delta G}{\delta h}$, the corresponding model

$$\frac{\partial h}{\partial t} = \nabla \cdot \left(M(\nabla h)\nabla e^{\frac{\delta G}{\delta h}}\right)$$

also needs further investigation.

**Acknowledgments.** This work was supported by the National Science Foundation under Grant No. DMS-1514826 and KI-Net RNMS11-07444.

*E-mail address*: maygao@ust.hk
*E-mail address*: hangjie@math.ucla.edu
*E-mail address*: jliu@math.duke.edu
*E-mail address*: witelski@math.duke.edu